\documentclass[a4paper,11pt]{article}

\usepackage{latexsym}
\usepackage{mathrsfs}
\usepackage{amssymb}
\usepackage{amsmath}
\usepackage{amsthm}
\usepackage[PostScript=dvips,nohug]{diagrams}
\usepackage[english]{babel}
\author{Yang Su}
\title{On the Classification of Certain Singular Hypersurfaces in $\cp$}
\date{}
\newtheorem{thm}{Theorem}
\newtheorem{lem}{Lemma}
\newtheorem{cor}{Corollary}
\newarrow{To}----{->}
\newarrow{Dashto}{}{dash}{}{dash}{->}
\newcommand{\CP}{(\mathbb{C}\mathrm{P}^{\infty})^{2}}
\newcommand{\cp}{\mathbb{C}\mathrm{P}^{4}}

\newcommand{\z}{\mathbb{Z} \oplus \mathbb{Z}}
\newcommand{\im}{\textrm{Im}}

\begin{document}
\maketitle

\section{Introduction}
Hypersurfaces in complex projective spaces defined by homogeneous polynomials  are important topological objects, arising naturally in
algebra, geometry and topology. It was first noted by Thom, the diffeomorphism type of smooth hypersurfaces depends only on the degree of the
defining polynomial; i.e., two $n$-dimensional smooth hypersurfaces in $\mathbb{C}\mathrm{P}^{n+1}$ are diffeomorphic if and only if they have equal degree. 
For singular hypersurfaces, the situation is more complicated. A theoretical solution to the classification problem of singular hypersurfaces is \cite{Di1}: 
let 
$P(n,d)$ be the moduli space of hypersurfaces of degree $d$ in $\mathbb{C}\mathrm{P}^{n}$, then there is a Whitney stratification on $P(n,d)$, such that 
two pairs $(\mathbb{C}\mathrm{P}^{n}, V_{f})$ and $(\mathbb{C}\mathrm{P}^{n}, V_{g})$ are topologically equivalent if the hypersurfaces $V_{f}$ and $V_{g}$ 
belong to the same connected component of a stratum of this stratification. 

Instead of considering the pair 
$(\mathbb{C}\mathrm{P}^{n}, V)$, in this paper we consider the classification of singular hypersurfaces as topological spaces. Since the nonsingular part of 
a hypersurface is a smooth manifold, it is natural to ask for a classification of singular hypersurfaces upto homeomorphisms, where the homeomorphisms are diffeomorphisms on the nonsingular parts. In this paper the hypersurfaces in $\cp$ with an isolated singularity are studied, and with some restrictions on the link of the singularity a classification in the above sense is obtained. The main result is the following:

\begin{thm}
For $i=1,2$, let $V_{i}\subset \mathbb{C}P^{4}$ be a hypersurface of degree $d_{i}$ with an isolated singularity $p_{i} \in V_{i}$, such that the link
of $p_{i}$ is diffeomorphic to $S^{2}\times S^{3}$. Suppose that the second integer homology group of the nonsingular part of  $V_{i}$ is isomorphic to
$\mathbb{Z} \oplus \mathbb{Z}$, and that $d_{i}$ is square-free. Let $\mu_{i}$ be the Milnor number of $p_{i}$. Then there is a
homeomorphism $f: V_{1} \to V_{2}$ which is a diffeomorphism on the nonsingular parts if and only if $d_{1}=d_{2}$ and $\mu_{1}=\mu_{2}$.
\end{thm}

As a consequence of this theorem, we have a classification of hypersurfaces with $A_k$ singularities.
\begin{cor}
For $i=1,2$, let $V_{i}\subset \mathbb{C}P^{4}$ be a hypersurface of degree $d_{i}$ with a unique singularity of type $A_{2k_{i}+1}$ ($k_{i} \ge 0$). If
$d_{i}<(k_{i}+5)/2$ and is square-free, then there is a homeomorphism $f: V_{1} \to V_{2}$ which is a diffeomorphism on the nonsingular parts if
and only if $d_{1}=d_{2}$ and $k_{1}=k_{2}$.

For $i=1,2$, let $V_{i}\subset \mathbb{C}P^{4}$ be a hypersurface of degree $d_{i}$ with a unique singularity of type $A_{2k_{i}}$. Then
there is a homeomorphism $f: V_{1} \to V_{2}$ which is a diffeomorphism on the nonsingular parts if and only if $d_{1}=d_{2}$ and $k_{1}=k_{2}$.
\end{cor}

It is interesting to compare this result with a known result on the topology of hypersurfaces with isolated singularities. As an application of the theoretical solution mentioned 
above, it is shown that if the degree is ``big enough'' compared with the number and the complexity of the singularities, then the space of the 
hypersurfaces with given degree
and singularities is connected and hence the topological type of these hypersurfaces is constant (\cite{Di1}). More precisely, let $d$ be the degree of $V$,
$k$ be the number of isolated singularities on $V$ and $s_{i},\dots, s_k$ be the $\mathscr{K}$-determinancy orders of the singularities, then the
space of such hypersurfaces is connected if $d \ge s_{1}+ \cdots + s_{k}+k-1$. The $\mathscr{K}$-determinancy order of an
$A_{k}$-singularity is $k+1$, hence if $d \ge 2k+2$, the pair $(d, k)$ is a complete invariant of hypersurface $V$ of degree $d$ with a unique
$A_{2k+1}$-singularity. On the other hand, corollary 1 gives information for $d$ relatively small ($d<(k+5)/2$).

The proof of theorem 1 consists of two parts: first we classify certain 6-manifolds with boundary in section 2; then in section 3 the invariants obtained in section 2 are computed for singular hypersurfaces. In section 3 the proof of corollary 1 is given, making use of a result on the cohomology of singular hypersurfaces in \cite{Di2}.  Examples of hypersurfaces fullfilling the assumptions will be given in the last section.

The results presented here are a part of the author's Ph.D. thesis at University Heidelberg. The author would like to thank Prof.~Matthias Kreck for his 
constant supervision and Dr.~Diarmuid Crowley for many useful discussions.

\section{Classification of certain 6-manifolds with boundary}
In this section we consider the classification of certain 6-manifolds with boundary. Let $M^{6}$ be a 6-dimensional oriented smooth manifold fulfilling the following
conditions:(A)
\begin{enumerate}
\item $M$ is simply-connected;
\item the boundary of $M$ is diffeomorphic to $S^{2}\times S^{3}$;
\item $H_2(M)$ is isomorphic to $\mathbb{Z} \oplus \mathbb{Z}$, $H_2(M, \partial M)$ is isomorphic to $\mathbb{Z}$ and the trilinear form
$$\begin{diagram}
H^2(M, \partial M)\times H^2(M, \partial M)\times H^2(M, \partial M) & \rTo & \mathbb{Z} \\
(x,y,z) & \mapsto & <x\cup y\cup z, [M, \partial M]> \\
\end{diagram}$$
is nontrivial.
\end{enumerate}

Our goal is to classify such manifolds upto orientation preserving diffeomorphisms.

Let us consider the invariants of such manifolds. First of all, there is a short exact sequence
$$0\to  H_2(\partial M)\to H_2(M)\to H_2(M, \partial M)\to 0.$$
Secondly, we have characteristic classes, namely, the Euler characteristic $\chi (M) \in \mathbb{Z}$, the Pontrjagin class $p_1(M) \in  H^{4}(M)$, and
the Stiefel-Whitney class $w_2(M)\in H^2(M;\mathbb{Z}/2)$. A generator $x$ of $H_2(M, \partial M)$ is called preferred
if the Kronecker dual of $x$, $x^*\in H^2(M, \partial M)$, satisfies that $\langle x^*\cup x^*\cup x^*, [M, \partial M] \rangle$ is positive.
Clearly the
preferred generator of $H_2(M, \partial M)$ is uniquely determined by this property and depends on the orientation. There is a symmetric bilinear form reflecting 
the cohomology multiplication of $H^*(M)$:
$$\begin{diagram}
q: & H^2(M)\times H^2(M) & \rTo & \mathbb{Z} \\
    & (u,v) & \mapsto & \langle u\cup v \cup x^{*}, [M, \partial M] \rangle .\\
\end{diagram}$$

Now we are able to formulate the classification of manifolds fulfilling condition (A) via
these invariants.

\begin{thm}
The Euler characteristic $\chi (M)$, the Poicar\'e dual of the Pontrjagin class $Dp_1(M)$, the Stiefel-Whitney class $w_2(M)$, the 
preferred generator $x\in H_2(M, \partial M)$, the short exact sequence
$$0\to  H_2(\partial M)\to H_2(M)\to H_2(M, \partial M)\to 0$$
and the bilinear form $q: H^2(M)\times H^2(M)\to \mathbb{Z}$ form a complete system
of invariants of oriented diffeomorphism type of $M$; i.e. there is an orientation preserving diffeomorphism between $M_0$ and $M_1$ if and only if
$\chi(M_0)=\chi(M_1)$ and there exists an isomorphism $\Phi$ between the short exact sequences
$$\begin{diagram}
0 & \rTo & H_2(\partial M_0) & \rTo & H_2(M_0) & \rTo & H_2(M, \partial M_0) & \rTo & 0\\
   &       & \dTo^{\Phi }          &        & \dTo^{\Phi } &    & \dTo^{\Phi }               &       & \\
0 & \rTo & H_2(\partial M_1) & \rTo & H_2(M_1) & \rTo & H_2(M, \partial M_1) & \rTo & 0\\
\end{diagram}$$
s.t. $\Phi (Dp_1(M_0))=Dp_1(M_1)$, $\Phi (x_0)=x_1$, and the dual of $\Phi $,
$\Phi ^{*}: H^2(M_1)\to H^2(M_0)$, is an isometry of the  bilinear forms $q_0$ on $H^2(M_0)$ and $q_1$ $H^2(M_1)$, and $\Phi ^{*}w_2(M_{1})=w_2(M_{0})$.
\end{thm}

We will use the modified surgery developed by M. Kreck in \cite{Kr} to prove this theorem. The modified surgery theory converts
a classification problem of manifolds into the problems of determining some bordism classes and certain obstruction. We will first deal with spin manifolds, 
and then show that with some trivial modification, the proof is also valid for nonspin manifolds.

\subsection{spin case}
Let $M$ be an oriented spin manifold, fulfilling the conditions (A). To apply the modified surgery, the first step here is to determine the normal 2-type of 
$M$.
Consider the fiberation
$$\xi: B=\CP \times B\textrm{Spin} \rTo^{\eta \times p} BO \times BO \rTo^{\oplus} BO,$$
where $p: B\textrm{Spin} \to BO$ is the canonical projection, $\eta:\CP \to BO$ is the classifying map of a trivial complex line bundle over $\CP$, and $\oplus $ is
the $H$-space structure of $BO$ given by the Whitney sum of universal vector bundles. Because $M$ is a simply-connected spin manifold, there is a unique 
classifying map $M\to B\textrm{Spin}$ of the spin structure on the stable normal bundle $\nu M$. By choosing an isomorphism 
$H_2(M) \stackrel{\sim}{\rTo} \mathbb{Z}\oplus \mathbb{Z}$
we get a map $M\to K(\mathbb{Z}\oplus \mathbb{Z}, 2)=\CP$,
which induces the given isomorphism.  Put these two maps together, we get a map
$$M\rTo^{\bar{\nu}} B\textrm{Spin}\times \CP,$$
which is clearly a lift of the normal Gauss map $\nu : M \to BO$: 
$$\begin{diagram}
& & B\\
& \ruTo^{\bar{\nu}} & \dTo^{\xi}\\
M & \rTo^{\nu} & BO\\
\end{diagram}$$

\noindent Since $\pi _2(B\textrm{Spin})=\pi _3(B\textrm{Spin})=0$, $\pi _3(\CP)=0$, we conclude that $(B, \xi)$ is the normal 2-type of $M$ and $\bar{\nu}$ is a normal 2-smoothing.

Now Let $M_i$ (i=0, 1) be as above, with the same Euler characteristic,
$\bar{\nu} _{i}:M_{i}\to B$ be a normal 2-smoothing of $M_i$ and
$$f:\partial M_{0}\to \partial M_1$$
be an orientation preserving diffeomorphism compatible with the normal 2-smoothings. Let $N_{f}=M_{0}\cup _{f}(-M_{1})$, then $\bar{\nu}_0$ and 
$\bar{\nu}_1$ fit together to give a map
$$\bar{\nu}_{f}=\bar{\nu}_{0} \cup \bar{\nu}_{1}:N_{f}\to B\mathrm{Spin}\times \CP.$$

The following lemma is a direct application of \cite[Corollary 4]{Kr}.

\begin{lem}
$f$ extends to an orientation preserving diffeomorphism
$$F:M_{0}\to M_{1}$$
compatible with the normal smoothings if and only if
$$[N_{f}, \bar{\nu}_{f}]=0\in \Omega_{6}(B;\xi)$$
\end{lem}

Here the bordism group $\Omega_{6}(B;\xi)$ is the bordism group of $6$-manifolds with normal $(B,\xi)$-structures. Since $\eta$ is the classifying map of a trivial bundle, the bordism group $\Omega_6(B;\xi)$ is identified with the spin bordism group $\Omega_6^{\mathrm{spin}}(\CP)$.

Therefore, it suffices to show that there exists an orientation preserving diffeomorphism
\mbox{$f:\partial M_{0}\to \partial M_1$} s.t. $(N_f, \bar{\nu}_{f})$ is null-bordant in
$\Omega_6^{\mathrm{spin}}(\CP)$. A standard computation with the Atiyah-Hirzebruch spectral sequence shows that there is an injective homomorphism
$$\begin{diagram}
\Omega_{6}^{\textrm{spin}}((\mathbb{C}P^{\infty})^2) & \rTo & H_2((\mathbb{C}P^{\infty})^2)\oplus H_6((\mathbb{C}P^{\infty})^2).\\
[Y,h] & \mapsto & (h_{*}Dp_{1}(Y), h_*[Y])\\
\end{diagram}$$

Now fix an orientation of $S^2\times S^3$ and choose an orientation reversing diffeomorphism $\varphi : \partial M_{0}\to S^2\times S^3$,
then the map
$$S^2\times S^3\rTo ^{\varphi ^{-1}} \partial M_{0} \subset M_0 \rTo^{\bar{\nu}_{0}} B \rTo^{\mathrm{pr}_{1}} \CP$$ 
extends to a map $S^2\times D^4\to \CP$ uniquely upto homotopy relative to the boundary. Let $Y_0=M_{0}\cup_{\varphi}(S^2\times D^4)$ and 
$$h_0:Y_0=M_{0}\cup_{\varphi}(S^2\times D^4) \to \CP$$
be the union of the corresponding maps on $M_0$ and $S^2\times D^4$. Since there is a unique spin structure on $S^{2}\times D^{4}$, we get an element $[Y_0, h_0]$ in $\Omega_{6}^{\mathrm{spin}}(\CP)$. (Here we choose the orientation of $S^2\times D^4$ so that it induces the fixed orientation of $S^2\times S^3$.) Do the same construction for the orientation reversing diffeomorphism $\partial M_1\rTo^{f^{-1}} \partial M_{0}\rTo^\varphi S^2\times S^3$
we obtain 
$$Y_1:=M_{1}\cup _{\varphi \circ f^{-1}}S^2\times D^4\rTo ^{h_1} \CP$$ 
and $[Y_1, h_1]\in \Omega_{6}^{\mathrm{spin}}(\CP)$. It is clear from the construction that
$$[N_{f}, \bar{\nu}_{f}]=[Y_0, h_0]-[Y_1, h_1]\in \Omega_{6}^{\mathrm{spin}}(\CP).$$

Now we study the bordism class $[Y_i, h_i]$. (For simplicity, we omit the subscription $i$ in the following discussion.) $Y$ is an oriented simply-connected 
6-manifold with $H_2(Y)\cong \mathbb{Z}\oplus \mathbb{Z}$ and the map $h:Y\to \CP$ induces an isomorphism on $H_2$. According to \cite{Wall}, the 
diffeomorphism type of $Y$ is determined by the Euler characteristic $\chi(Y)$, (here $\chi(Y)=\chi(M)+2$),
the Poincar\'e dual of the Pontrjagin class $Dp_{1}(Y)$ and the trilinear form
$$\begin{diagram}
H^2(Y)\times H^2(Y)\times H^2(Y) &\rTo^{\ \ \ \ \ \mu} & \mathbb{Z}.\\
(a, b, c)                                         & \mapsto     & \langle a\cup b\cup c, [Y]\rangle \\
\end{diagram}$$

It is seen from the Mayer-Vietoris sequence that the inclusion
$j:M\to Y$
induces an isomorphism on $H_2$. We identify $H_2(M)$ and $H_2(Y)$ using this isomorphism.
Choose a basis of $H_2(M)$, $\{e_1, e_2\}$, s.t. $e_1$ is the image of a generator of $H_2(\partial M)$ under the inclusion $H_2(\partial M)\to H_2(M)$,
and $e_2$ maps to the preferred generator of $H_2(M, \partial M)$ under the projection $H_2(M)\to H_2(M, \partial M)$:

$$\begin{diagram}
H_2(\partial M) & \rTo & H_2(M) & \rTo & H_2(M, \partial M)\\
1              &\mapsto & e_1,\ e_2 & \mapsto & x\\ 
\end{diagram}$$

\noindent Under this basis, the invariants of $Y$ can be expressed as follows ($\star$):

\begin{itemize}
\item[-] $\chi (Y)=\chi (M)+2$
\item[-] $Dp_{1}(Y)=p\cdot e_1+Dp_{1}(M)\cdot e_2$, for some  $p\in \mathbb{Z}$
\item[-] $\mu (e_1^*, e_1^*, e_1^*)=\lambda$, for some $\lambda \in \mathbb{Z}$
\item[-] the restriction of $\mu$ on $H^2(Y)\times H^2(Y)\times \mathbb{Z}\cdot e_{2}^{*}$
\end{itemize}

\noindent Here $Dp_{1}(M)$ is understood as an integer under the isomorphism $H_2 (M, \partial M)\\ \cong \mathbb{Z}$ given by the preferred generator,
$e_{i}^{*} \in H^2(Y)$ is the Kronecker dual of $e_{i}$. According to \cite{Wall}, 
$p\equiv 4\lambda \pmod{24}$. The restriction of $\mu$ on $H^2(Y)\times H^2(Y)\times \mathbb{Z}\cdot e_{2}^{*}$ is equivalent to the bilinear form $q$.

Concerning the relation between these invariants and the bordism class $[Y, h]$, we have the following lemma.

\begin{lem}\label{pt}
Let $Y^6$ be a spin $6$-manifold, $h:Y\to \CP$ be a map inducing an isomorphism on $H_2$. Then the bordiam class 
$[Y, h]\in \Omega_{6}^{\mathrm{spin}}(\CP)$ is determined by the Poincar\'e dual of the Pontrjagin class $Dp_{1}(Y)$ and the trilinear
form
$$\begin{diagram}
H^2(Y)\times H^2(Y)\times H^2(Y) & \rTo^{\ \ \ \ \ \mu} & \mathbb{Z}.\\
(a, b, c)                        & \mapsto     & \langle a\cup b\cup c, [Y]\rangle \\
\end{diagram}$$
\end{lem}

\begin{proof}
The bordism class $[Y,h]$ is determined by $Dp_1(Y)$ and $h_*[Y]$. Therefore we only need to show that the trilinear form $\mu$ determines $h_*[Y]\in H_6(\CP)$. Since $H^6(\CP)$ is generated by elements of the form
$u\cup v\cup w$, $u, v, w\in H^2(\CP)$, $h_*[Y]$ is determined by the evaluation
$$\langle u\cup v\cup w, h_*[Y]\rangle =\langle h^*(u)\cup h^*(v)\cup h^*(w), [Y]\rangle.$$

\noindent $h_*$ is an isomorphism, so is $h^*$, hence the evaluation is equivalent to the trilinear form $\mu$.
\end{proof}

If we identify $H^2(M)$ and $H^2(Y)$ via $j^*$, then $\mu$ can be viewed as a trilinear form on $H^2(M)$. A big part of this trilinear form is determined by the bilinear form $q$ on $H^2(M)$. More precisely, we have

\begin{lem}
The restriction of $\mu$ on $H^2(M)\times H^2(M)\times H^2(M, \partial M)$ is equivalent to the bilinear form $q:H^2(M)\times H^2(M)\to \mathbb{Z}$.
\end{lem}

\begin{proof}
We have the following commutative diagram
$$\begin{diagram}
H^2(Y) & \times & H^2(Y) & \times & H^2(Y) & \rTo^\mu & \mathbb{Z}\\
\uTo^\cong &     & \uTo^\cong &     &  \uTo^\cong &        &  \Vert \\
H^2(M) & \times & H^2(M) & \times & H^2(M) & \rTo^\mu  & \mathbb{Z}\\
\Vert      &          &  \Vert    &          & \bigcup  &              & \Vert \\
H^2(M) & \times & H^2(M) & \times & H^2(M, \partial M) & \rTo & \mathbb{Z}
\end{diagram}$$
The map on the bottom is $(u, v, w)\mapsto \langle u\cup v\cup w, [M, \partial M]\rangle$, for $u, v\in H^2(M)$, $w\in H^2(M, \partial M)$. This is equivalent to
the bilinear form $q$.
\end{proof}

\begin{lem}\label{change}
Let $Y=M\cup_{\varphi}(S^2\times D^4)$ and $h$ be as above. Then there exists an orientation preserving diffeomorphism $g:S^2\times S^3\to S^2\times S^3$ s.t. for $Y'=M\cup_{g\circ \varphi}S^{2} \times D^{4}$, we have $p=\lambda=0$ in $(\star)$.
\end{lem}

\begin{proof}
Let $P=(-(S^2\times D^4))\cup_g (S^2\times D^4)$ for some orientation preserving diffeomorphism $g:S^2\times S^3\to S^2\times S^3$, then $P$ is
a 6-dimensional, simply-connected spin manifold with $H_2(P)\cong \mathbb{Z}$. Define a map $k:P\to \CP$ as follows:

On the first copy of $S^2\times D^4$, $k$ is the extension of
$$S^2\times S^3\rTo^{\varphi^{-1}} \partial M \subset M\rTo^{\bar{\nu}} B \rTo^{pr_{1}} \CP,$$
and on the second copy of $S^2\times D^4$ $k$ is the extension of
$$S^2\times S^3\rTo^g S^2\times S^3\rTo^{\varphi^{-1}} \partial M \subset M \rTo^{\bar{\nu}} B \rTo^{pr_{1}} \CP.$$

\noindent Then from the construction, it is seen that $$k_*:H_2(P)\rTo^\cong h_*(\mathbb{Z}\cdot e_{1}) \subset H_2(\CP)$$ and
$$[Y, h]+[P, k]=[Y', h']\in \Omega_{6}^{\mathrm{spin}}(\CP).$$

\noindent Therefore, $h_{*}Dp_{1}(Y)+k_{*}Dp_{1}(P)=h_{*}'Dp_{1}(Y')$ and
$$\begin{array}{ll}
 & \mu' (e_1^*, e_1^*, e_1^*)\\
  = & \langle e_1^*\cup e_1^*\cup e_1^*, [Y']\rangle\\
   = &  \langle (h'^*)^{-1} (e_1^*\cup e_1^*\cup e_1^*), h'_{*}[Y']\rangle\\
   = &  \langle(h^{*})^{-1} (e_1^*\cup e_1^*\cup e_1^*), h_*[Y]\rangle + \langle (h^{*})^{-1} (e_1^*\cup e_1^*\cup e_1^*), k_*[P]\rangle \\
   = & \mu (e_1^*, e_1^*, e_1^*)+\langle u^*\cup u^*\cup u^*, [P]\rangle
\end{array}$$

\noindent where $u^* =k^* ((h^* )^{-1}(e_1 ^*))\in H^2 (P)$ is a generator.

According to the classification result of \cite{Wall}, for any $b, c\in \mathbb{Z} $, there exists an orientation preserving
diffeomorphism $g:S^2\times S^3\to S^2\times S^3$ s.t. $P=(-(S^2\times D^4))\cup_g (S^2\times D^4)$ satisfies

$$\left \{ \begin{array}{ccl}
Dp_{1}(P) & = & 4b\cdot u\\
\\
\langle u^*\cup u^*\cup u^*, [P]\rangle & = & 6c+b
\end{array} \right.$$

\noindent where $u\in H_2(P)$ is the Kronecker dual of $u^*$. By choosing
$$\left \{ \begin{array}{rcl}
b & = & -p/4\\
\\
c & = & (p-4\lambda)/24\\
\end{array} \right.$$

\noindent we get a corresponding $g$. For the corresponding $Y'$ we have
$$\left \{ \begin{array}{ccl}
Dp_{1}(Y') & = & p\cdot e_1+Dp_{1}(M)\cdot e_2-p\cdot e_1 =  Dp_{1}(M)\cdot e_2\\
\\
\mu'(e_1^*, e_1^*, e_1^*) & = & 0
\end{array} \right.$$

\noindent This proves the lemma.
\end{proof}

\begin{proof}[Proof of Theorem 2]
Suppose that $M_0$ and $M_1$ satisfy the assumptions. We choose normal 2-smoothings $\bar{\nu}_{i} \ (i=1,2)$, compatible with $\Phi$, 
i.e.~the isomorphism between $H_2(M_0)$ and $H_2(M_1)$ induced by $\bar{\nu}_0$ and $\bar{\nu}_1$ coinsides with $\Phi$. 
Lemma \ref{change} ensures us to choose diffeomorphisms \mbox{$\varphi_{i}:\partial M_i \to S^2\times S^3$} $(i=1,2)$, s.t. $Y_{i}=M_{i}\cup_{\varphi_{i}} (S^2\times D^4)$
satisfies
$$\left \{ \begin{array}{ccl}
Dp_{1}(Y_{i}) & = & Dp_{1} (M_{i} )\cdot e_{2}^{i}\\
\\
\mu (e_{1}^{i*}, e_{1}^{i*}, e_{1}^{i*}) & = & 0\\
\end{array}\right.$$

\noindent where $\{ e_{1}^{0},\ e_{2}^{0}\}$ is a basis of $H_2(M_0)$ as before and $\{ e_{1}^{1}=\Phi (e_{1}^{0}),\ e_{2}^{1}=\Phi (e_{2}^{0})\}$
is a basis of $H_2(M_1)$. Let $f$ be the composition
$$\partial M_0\rTo^{\varphi_0} S^2\times S^3\rTo^{\varphi_{1}^{-1}} \partial M_1,$$

\noindent then the isomorphism between $H_2(\partial M_0)$ and $H_2(\partial M_1)$ induced by $f$ coinsides with $\Phi$. Thus $f$ is compatible with the 
normal 2-smoothings since $\Phi$ is an isomorphism of the short exact sequences. We claim that
$[Y_0 , h_0 ]=[Y_1 , h_1 ]\in \Omega_6 ^{\mathrm{spin}}(\CP)$. Once we have shown this, since
$[N_{f}, \bar{\nu_f}]=[Y_0, h_0]-[Y_1, h_1]$, lemma 1 implies that $M_0$ is diffeomorphic to $M_1$. The diffeomorphism induces the same isomorphism 
on homology as $\Phi$ does, hence is compatible with the $B$-structures. 

Since $\Phi :H_2(M_0 , \partial M_0 )\to H_2(M_1 , \partial M_1 )$ maps $Dp_{1}(M_0 )$ to $Dp_{1}(M_1 )$ and $\Phi ^* :H^2 (M_1 )\to H^2 (M_0 )$ is an isometry of the bilinear forms, it follows that
$\Phi :H_2 (Y_0 )\to H_2 (Y_1 )$ maps $Dp_{1} (Y_0 )$ to $Dp_{1} (Y_1 )$ and \mbox{$\Phi ^* :H^2 (Y_1 )\to H^2 (Y_0 )$} preserves the trilinear form. 
By lemma \ref{pt}, $[Y_0 , h_0 ]=[Y_1 , h_1 ]\in \Omega_6 ^{\mathrm{spin}}(\CP)$. This finishes the proof of theorem 2 for spin manifolds.
\end{proof}

\subsection{nonspin case}
The proof of theorem 2 for nonspin manifolds is essentially the same as that for spin manifolds, except for some minor modification. 

First of all, in this case, the normal $2$-type of $M$ is described as follows: Consider the fibration
$$\xi: B=\CP \times B\mathrm{Spin} \rTo^{\eta \times p} BO \times BO \rTo^{\oplus} BO,$$
where $p: B\mathrm{Spin} \to BO$ is the canonical projection, $\eta:\CP \to BO$ is the classifying map of the complex line bundle $\mathrm{pr}_{1}^{*}(H)$, where
$H$ is the canonical line bundle over $\mathbb{C}\mathrm{P}^{\infty}$ and $\mathrm{pr}_{1}: \CP \to \mathbb{C}\mathrm{P}^{\infty}$ is the projection
to the first factor. Let $\bar{\nu}_{1}: M \to \CP$ be a map which induces an isomorphism on $H_{2}$ and s.t.~$\bar{\nu}_{1}^{*}((1,0)) \equiv w_{2}(\nu M) \pmod {2}$. (This is the case since $M$ is nonspin.)
Then $w_{2}(\nu M-\bar{\nu}_{1}^{*}(\eta))=0$ and there is a (unique) lift $\bar{\nu}_{2}: M \to B\mathrm{Spin}$ classifying the stable bundle 
$\nu M-\bar{\nu}_{1}^{*}(\eta)$.
Let $\bar{\nu}=\bar{\nu}_{1} \times \bar{\nu}_{2}: M \to B$, then $\bar{\nu}$ is a normal $2$-smoothing of $M$ in the normal $2$-type $(B,\xi)$.

The construction of the manifold $Y=M \cup_{\varphi} (S^{2} \times D^{4})$ and the map $h$ gives an element in the corresponding bordism group.
Lemma 1 holds for this normal structure as well. Now the bordism group $\Omega_{6}(B;\xi)$ is identified with the twisted spin bordism group 
$\Omega_{6}^{\textrm{spin}}(\CP ; \eta)$, which is the bordism group of maps
$f$ from closed $6$-manifolds $X$ to $\CP$, together with a spin structure on $f^{*}(\eta) \oplus \nu X$. There is an isomorphism
$$\Omega_{6}^{\textrm{spin}}(\CP ; \eta) \cong \widetilde{\Omega}_{8}^{\textrm{spin}}(\textrm{Th}(\eta)),$$
where $\textrm{Th}(\eta)$ is the Thom space of the complex line bundle $\eta$. A computation with the Atiyah-Hirzebruch spectral sequence shows that there is an injective homomorphism
$$\begin{diagram}
\Omega_{6}^{\textrm{spin}}(\CP ; \eta) & \rTo    & H_{2}(\CP) \oplus H_{6}(\CP). \\
[Y,h]                         & \mapsto & (h_{*}Dp_{1}(\tau Y \oplus h^{*}(\eta)), h_{*}[Y])\\
\end{diagram}$$
Since
$$\begin{array}{rcl}
h_{*}Dp_{1}(\tau Y \oplus h^{*}(\eta)) & = & h_{*}Dp_{1}(Y) + h_{*}Dp_{1}(h^{*}\eta)\\
                                                     & = & h_{*}Dp_{1}(Y) + h_{*}Dh^{*}p_{1}(\eta)\\
                                                     & = & h_{*}Dp_{1}(Y) + p_{1}(\eta)\cap h_{*}([Y]),\\
\end{array}$$
the bordism class is determined by $h_{*}Dp_{1}(Y)$ and $h_{*}([Y])$.

Finally note that there is still the relation $p \equiv 4\lambda  \pmod {24}$ and lemma 4 is also valid in this case.
One easily checks that after the modification the proof for spin manifolds is valid for  nonspin manifolds.

\section{Computing the invariants for singular hypersurfaces}
In last section a complete system of diffeomorphism invariants for the $6$-manifolds under consideration is obtained. In this section we will compute these invariants for
the smooth part of hypersurfaces in $\cp$. Then theorem 1 follows from theorem 2 and this computation.

Let $V \subset \cp$ be a hypersurface of degree $d$, with a unique singularity $p$. Because of the conic structure of the hypersurface near the
singularity $p$, there exits a small open ball $D_\epsilon \subset \cp$ with center at $p$, such that $\overline{D}_\epsilon \bigcap V$ is homeomorphic to the cone over
$\partial \overline{D}_\epsilon \bigcap V$.
$\partial \overline{D}_\epsilon \bigcap V$ is called the link of the singularity $p$ (c.f.~\cite{Di1}).
Let $M = V-B$, then $M$ is a smooth manifold with boundary $\partial \overline{D}_\epsilon \bigcap V$, the interior of $M$ is diffeomorphic to the
nonsingular part of $V$, and $V$ is homeomorphic to $M/ \partial M$. We call $M$ the smooth part of $V$ and have the following

\begin{lem}
Let $M$ be the smooth part of $V$, then the Euler characteristic  $\chi(M)$ and the first Pontrjagin class  $p_{1}(M)$ are
determined by the degree $d$ and the Milnor number $\mu_{p}$ of the singularity $p$. $M$ is spin if and only if $d$ is odd.
\end{lem}

\begin{proof}
At first consider a smooth hypersurface $i: V_{0}\subset \cp $ of degree $d$. Let $\nu(i)$ be the normal bundle of the embedding, then
$TV_{0}\oplus \nu(i) =i^{*}T\cp$. It is known that $\nu(i) =i^{*}(H^{\otimes d})$, where $H$ is the canonical complex line bundle over $\cp$. Thus
$c(V_{0})=i^{*}c(T\cp)/i^{*}c(H^{\otimes d})=i^{*}((1+x)^{5}/(1+dx))$, where $x=c_{1}(H)$. Therefore
$\chi(V_{0})=\langle c_{3}(TV_{0}), [V_{0}]\rangle$, $p_{1}(V_{0})=-c_{2}(TV_{0}\oplus \overline{TV}_0)$ are determined by $d$, where $\overline{TV}_0$ is the conjugate bundle of $TV_0$.
Since $c_1(V_0)=(5-d)i^*(x)$ and $c_{1}(V_{0})\equiv w_{2}(V_{0}) \ (\textrm{mod}\ 2)$, $V_{0}$ is spin if and only if $d$ is odd. 

Now let $V_{0}$ be a small deformation of $V$ such that $V_{0}$ is smooth.
Then $\overline{D}_\varepsilon\cap V_{0}$ can be identified
with the closed Milnor fiber of the singularity and $M$ is diffeomorphic to $V_{0}-D_\varepsilon$ (c.f.~\cite[page 163]{Di1}). Then
$\chi(M)=\chi(V_{0}) - \chi(V_{0} \cap \overline{D}_\varepsilon)$. Since the Milnor fiber is homotopy equivalent to a wedge of $\mu$ copies of $3$-spheres, the Euler characteristic of $V_{0} \cap \overline{D}_\varepsilon$ is $1-\mu_{p}$. Therefore $\chi(M)$ is determined by $d$ and
$\mu_{p}$. 

If we identify $M$ and $V_{0}-D_\varepsilon$, then the inclusion $j: M \to V_{0}$ induces an isomorphism $j^{*}: H^{4}(V_{0}) \to H^{4}(M)$, therefore 
$p_{1}(M) =j^{*}p_{1}(V_{0})\in H^{4}(M)\cong \mathbb{Z}$ is determined by $d$. 

If $d$ is odd, then $V_{0}$ is spin and therefore $M$ is spin.
If $d$ is even, $V_{0}$ is nonspin. We show that $M$ is nonspin: if $M$ is spin, since $\overline{D}_\varepsilon\cap V_{0}$ is also spin and there is a unique spin structure on $\partial M$ ($\partial M$ is simply-connected), the spin
structures on $M$ and $\overline{D}_\varepsilon\cap V_{0}$ fit together to give rise to a spin structure on $V_{0}$, which is a contradiction. Therefore $M$ is nonspin.
\end{proof}

Next we compute the homology and cohomology of the smooth part $M$. By Lefschetz theorem, $V$ is simply-connected and $H_2(V)$ is isomorphic to $\mathbb Z$, therefore by Van-Kampen theorem $M$ is simply-connected and there is an exact sequence
$$H_2(\partial M)\rTo^{i_*} H_2(M)\rTo H_2(M, \partial M)\rTo 0.$$ 

\noindent Since $H_2(M, \partial M)\cong H_2(V)\cong\mathbb Z$, $H_2(M)$ is isomorphic to $\mathbb Z\oplus \im i_*$. From now on we assume that 
$\partial M$ is diffeomorphic to $S^{2}\times S^{3}$ and $\im i_*\cong \mathbb Z$. This assumption is fulfilled if $p$ is an $A_{2k+1}$-singularity and $d<(k+5)/2$.

Let $V^{*}= V-\{p\}$ be the nonsingular part and
$U= \cp - V$ be the complement. Since the embedding $i: V^{*}\hookrightarrow \cp -\{p\}$ is proper, there is a Gysin sequence (c.f.~\cite[page 314, 321]{Do}):
$$\cdots \to H^{k}(\cp -\{p\}) \stackrel{j^*}{\to} H^{k}(U)\stackrel{R}{\to} H^{k-1}(V^{*})\stackrel{\delta}{\to} H^{k+1}(\cp -\{p\})\to \cdots,$$

\noindent where $j: U \subset \cp -\{p\}$ denotes the inclusion, and the homomorphism $R$ is the so-called Poincar\'e-Lerray residue \cite{Di1}. For $k=3$, we have a commutative diagram:
$$\begin{diagram}
0 & \rTo & H^{3}(U) & \rTo^{R} & H^{2}(V^{*}) & \rTo^{\delta} & H^{4}(\cp -\{p\})\cong \mathbb{Z} \\
   &        &               &               &                    & \luTo^{i^{*}} & \uDashto                          \\
   &       &                &              &                       &                   & H^{2}(\cp -\{p\}) \cong \mathbb{Z} \\
\end{diagram}$$

\noindent The composition $\delta \circ i^{*}: H^{2}(\cp -\{p\}) \to H^{4}(\cp -\{p\})$ is a multiplication by $d$, and by the assumption above, $H^2(V^*)\cong H^2(M) \cong \mathbb{Z} \oplus \mathbb{Z}$, therefore $H^{3}(U)$ is isomorphic to $\mathbb{Z}$. Let $u \in H^{3}(U)$ be a generator. We have the following

\begin{lem}\label{cup}
For any $y \in H^{2}(V^{*})$, the cup product $y \cup R(u)$ is  $0$.
\end{lem}

\begin{proof}
Let $T$ be a tubular neighbourhood of $V^{*}$ in $\cp -\{p\}$, $T_{0}=T-V^{*}$ be the complement of the zero section and $j_{0}: T_{0} \to T$ be the inclusion. Since $H^2(V^*)$ and $H^4(V^*)$ are torsion free, it suffices to compute $y \cup R(u)$ with $\mathbb{Q}$-coefficients.
We consider the associated Gysin sequence with $\mathbb{Q}$-coefficients:
$$\cdots \rTo H^{k}(T_{0}) \rTo^{R_{0}} H^{k-1}(T) \rTo H^{k+1}(T) \rTo^{j_{0}^{*}} H^{k+1}(T_{0}) \rTo \cdots,$$
where $R_{0}$ is the composition
$$H^{k}(T_{0}) \rTo^{\delta} H^{k+1}(T, T_{0}) \rTo^{\varphi} H^{k-1}(T),$$
where $\varphi$ is the Thom isomorphism.

For $k=3$, the map $H^{2}(V^{*})\cong H^{2}(T) \to H^{4}(T)\cong H^{4}(V^{*})$ is just the cup product with the Euler class of the normal bundle,  hence
 is surjective.

Let $b \in H^{1}(T_{0}) \cong \mathbb{Z}$ be a generator, then for any $x \in H^{3}(T_{0})$, we have
$$\begin{array}{rcl}
R_{0}((j_{0}^{*}R_{0}(x))\cup b) & = &\varphi \circ \delta ((j_{0}^{*}R_{0}(x))\cup b) \\
                                 &  = & \varphi (R_{0}(x)\cup \delta(b)) \\
                                 &  = & R_{0}(x) \cup R_{0}(b) =R_{0}(x),\\
\end{array}$$

\noindent since $\delta$ is an $H^{*}(T_{0})$-mod map and the Thom isomorphism $\varphi$ is an  $H^{*}(T)$-mod map.

Thus $j_{0}^{*}R_{0}(x)\cup b-x=j_{0}^{*}(z)$ for some $z \in H^{3}(T)$. For any $y \in H^{2}(T)$, we claim that $y \cup R_{0}(x)=0$. This
follows from the following calculation:
$$\begin{array}{rcl}
y \cup R_{0}(x) & = &\varphi (y \cup \delta(x))\\
                       & = & \varphi \delta (j_{0}^{*}(y) \cup x) \\
                       & = & R_{0}(j_{0}^{*}(y) \cup x)\\
                       & = & R_{0}(j_{0}^{*}(y) \cup (j_{0}^{*}R_{0}(x)\cup b - j_{0}^{*}(z))) \\
                       & = & R_{0}(j_{0}^{*}(y \cup R_{0}(x))\cup b -j_{0}^{*}(y \cup z)).\\
\end{array}$$

\noindent Note that $j_{0}^{*}(y \cup R_{0}(x))=0$ since $j_{0}^{*}: H^{4}(T) \to H^{4}(T_{0})$ is a trivial map and $y \cup z=0$ since $H^{5}(T)=0$. Therefore $y \cup R_{0}(x)=0$.

The statement of this lemma is then proved by the commutative diagram
$$\begin{diagram}
H^{3}(U)   & \rTo ^{R}     & H^{2}(V^{*}) \\
\dTo       &               & \uTo^{\cong}       \\
H^{3}(T_{0}) &\rTo^{R_0}   & H^2(T)        \\
\end{diagram}$$

\noindent where the homomorphism $H^3(U) \to H^3(T_0)$ is induced by the inclusion \mbox{$T_0 \subset U$}.
\end{proof}

Now we are able to prove theorem 1.

\begin{proof}[Proof of Theorem 1]
Let $V \subset \cp$ be a hypersurface fulfilling the assumptions. Let $M$ be the smooth part of $V$. Then $M$ fulfills the conditions
(A). Consider the short exact sequence
$$0\rTo H_{2}(\partial M) \rTo^{i} H_{2}(M) \rTo^{j} H_{2}(M, \partial M)\rTo 0.$$

\noindent Let $y\in H_{2}(\partial M) \cong \mathbb{Z}$ be a generator and $x\in H_{2}(M, \partial M) \cong \mathbb{Z}$ be the preferred generator.
Let $b=i(y)$ and $a\in j^{-1}(x)$, then $\{a,b\}$ is a basis of $H_{2}(M)\cong \z$. Let $a^{*}$, $b^{*}$ be the Kronecker dual of $a$, $b$, then 
$\{a^{*}, b^{*}\}$ is a basis of $H^{2}(M)$. Recall we have a Gysin sequence 
$$0 \rTo H^3(U) \rTo^R H^2(V^*) \rTo^{\delta} H^4(\mathbb{C}\mathrm{P}^4-p)$$
and $H^3(U) \cong \mathbb{Z}$. Let $u \in H^3(U)$ be a 
generator, under the identification of $H^{2}(V^{*})$ and $H^{2}(M)$ by the inclusion, the primitive element $R(u)$ can be written as $R(u)=ma^{*}+nb^{*}$
for some coprime $m,n\in \mathbb{Z}$. According to lemma \ref{cup}, and the fact that $q(a^*,a^{*})=a^*\cup a^{*}\cup x^*=d$, we
have
$$\left \{ \begin{array}{rcccl}
0 & = & q(R(u),a^{*}) & = & md+n \cdot q(b^{*},a^{*}) \\
\\
0 & = & q(R(u),b^{*}) & = & m \cdot q(a^{*},b^{*})+n \cdot q(b^{*},b^{*})\\
\end{array}\right.$$

\noindent This implies
$$\left \{ \begin{array}{rcl}
q(a^{*},b^{*}) & = & -md/n\in \mathbb{Z} \\
\\
q(b^{*}, b^{*}) & = & m^2d/(n^2)\in \mathbb{Z}\\
\end{array}\right.$$

\noindent Since $m$ and $n$ are coprime, $d$ is divisible by $n^{2}$. If $d$ is square-free, then $n$ equals to $1$ and thus $R(u)=ma^{*}+b^{*}$.
We can perform a basis change in $H_{2}(M)$ by replacing $a$ by $a+mb$. Under this new basis we have $R(u)=b^{*}$ and the bilinear form $q$ is represented
by the matrix
$$\left( \begin{array}{cc}
d & 0\\
0 & 0
\end{array} \right).$$

Furthermore, by lemma 5, $p_1(M)$ and $\chi(M)$ are determined by the degree $d$ and the Milnor number $\mu$ of $p$. $M$ is spin if and only if $d$ is odd. 
Note that $i^{*}(w_{2}(M))=0$ since $\partial M$ is spin. Thus if $d$ is even, $w_{2}(M)$ is the image of the generator of $H^{2}(M, \partial M; \mathbb{Z}/2)$.
Let $V_{1}$ and $V_{2}$ be hypersurfaces in theorem 1, with degree $d_{i}$ and Milnor number $\mu_{i}$. Let $M_{1}$
(resp. $M_{2}$) be the smooth part of $V_{1}$ (resp. $V_{2}$). Let $\{ a_{1}, b_{1} \}$ (resp. $\{ a_{2}, b_{2} \}$) be the basis of $H_{2}(M_{1})$
(resp. $H_{2}(M_{2})$) chosen as above. If $d_{1}=d_{2}$ and $\mu_{1}=\mu_{2}$, then $\chi(M_{1}) = \chi(M_{2})$. Let $\Phi$ be an isomorphism of the 
corresponding short exact sequences, mapping $a_{1}$ to
$a_{2}$ and $b_{1}$ to $b_{2}$. Thus $Dp_{1}$, $w_{2}$ and $q$ are all  preserved by $\Phi^{*}$. Therefore according to theorem 1 $M_{1}$ is diffeomorphic to $M_{2}$.
This diffeomorphism extends to a homeomorphism between $V_{1}$ and $V_{2}$. Conversely, if there is a homeomorphism between $V_{1}$ and $V_{2}$ which is
diffeomorphism on the nonsingular parts, then $M_{1}$ and $M_{2}$ are diffeomorphic, thus $d_{1}=d_{2}$ and $\mu_{1}=\mu_{2}$.
\end{proof}

\begin{lem}
The link of a singularity  $(X,0)\subset (\mathbb{C}^{4},0)$ of type $A_{k}$ is $S^{2}\times S^{3}$ for $k$ odd and $S^{5}$ for $k$ even.
\end{lem}
\begin{proof}
A surface singularity in $\mathbb{C}^{3}$ of type $A_{k}$ is defined by the equation
$$x_{1}^{2}+x_{2}^{2}+x_{3}^{k+1}=0.$$
The resolution graph of this singularity is the same as the Dynkin diagram of the Lie algebra $A_{k}$. A singularity of type $A_{k}$
in $\mathbb{C}^{4}$ defiend by the equation 
$$x_{1}^{2}+x_{2}^{2}+x_{3}^{2}+x_{4}^{k+1}=0$$
is the stabilization of a surface $A_{k}$-singularity. Therefore under a distinguished basis,
the Milnor lattice of this singularity is represented by the matrix
$$\left( \begin{array}{ccccc}
0      &   1    &      0 & \ldots & \ldots \\
-1     &   0    &      1 & \ldots & \ldots \\
0      &  -1    &      0 & \ldots & \ldots \\
\vdots & \vdots & \vdots & \ddots &   1 \\
\vdots & \vdots & \vdots &    -1  &     0
\end{array} \right)
$$
After a simple calculation it is seen that this bilinear form is equivalent to
$$\bigoplus_{k/2} \left( \begin{array}{cc}
0 & 1\\
-1 & 0 \\
\end{array} \right )
$$
for $k$ even and 
$$\bigoplus_{(k-1)/2} \left( \begin{array}{cc}
0 & 1\\
-1 & 0 \\
\end{array} \right ) \bigoplus (\mathbb{Z}, ()_{0})$$
for $k$ odd, where $(\mathbb{Z}, ()_{0})$ denotes the trivial form on $\mathbb{Z}$. Then by the classification of
simply-connected $5$-manifolds of Smale, we see that the link of an $A_{k}$-singularity is  $S^{2}\times S^{3}$ for $k$ odd and $S^{5}$ for $k$ even.
\end{proof}

\begin{lem}
Let $V \subset \cp$ be a hypersurface of degree $d$ with a unique singularity $p$ of type $A_{2k+1}$, $M$ be the smooth part. Then if $d<(k+5)/2$,
$H_{2}(M)\cong \mathbb{Z} \oplus \mathbb{Z}$; if $d \ge (k+5)/2$, the rank of $H_{2}(M)$ is $1$.
\end{lem}
\begin{proof}
Since there is an exact sequence
$$H_2(\partial M)\to H_2(M)\to H_2(M, \partial M)\to 0$$
and  $H_2(\partial M) \cong H_2(M, \partial M) \cong \mathbb{Z}$, to prove $H_{2}(M)\cong \mathbb{Z} \oplus \mathbb{Z}$, it suffices to show that
$\dim_{\mathbb{C}} H_{2}(M;\mathbb{C})=2$. Since the local equation of an $A_{2k+1}$-singularity is
$$x_{1}^{2}+x_{2}^{2}+x_{3}^{2}+x_{4}^{2k+1}=0,$$
there exists a hyperplane $H$ in $\cp$, intersecting $V$ transversely, such that $p \in H$ and $(V \bigcap H, p)$ is a singularity of type $A_{1}$.
Without lose of generality, we may suppose $H$ is defined by the equation $x_{0}=0$. Consider the linear system
$$\mathscr{S}_{k}=\{h \in \bar{S}_{2d-k-5}| h(q)=0 \},$$
where $\bar{S}$ denotes the polynomial ring $\mathbb{C}[x_{1}, \cdots , x_{4}]$ and $\bar{S}_{m}$ the homogeneous part of degree $m$. Then
according to \cite[proposition 3.4]{Di2}, 
$$\dim_{\mathbb{C}} H_{2}(M;\mathbb{C})=\dim_{\mathbb{C}} H^{4}(V;\mathbb{C})=1+ \textrm{def} \mathscr{S}_{k}=1+1- \textrm{codim} \mathscr{S}_{k}.$$
It is clear that $\textrm{codim} \mathscr{S}_{k}=0$ if $d<(k+5)/2$ and $\textrm{codim} \mathscr{S}_{k}=1$ if $d \ge (k+5)/2$.
\end{proof}

\

\begin{proof}[Proof of Corollary 1]
Let $V\subset \mathbb{C}P^{4}$ be a hypersurface of degree $d$ with a unique singularity of type $A_{2k+1}$ ($k \ge 0$), then the Milnor number $\mu$ of this 
singularity is $-k+1$ and according to lemma 7, the boundary of the smooth part $M$ is diffeomorphic to $S^{2} \times S^{3}$.
If $d<(k+5)/2$, lemma 8 implies $H_{2}(M) \cong \mathbb{Z} \oplus \mathbb{Z}$. Then corollary 1 follows from theorem 1.

If $V\subset \mathbb{C}P^{4}$ is a hypersurface
of degree $d$ with a unique singularity of type $A_{2k}$, then the boundary of the smooth part $M$ is diffeomorphic to $S^{5}$. In this case
$H_{2}(M) \cong H_{2}(M, \partial M) \cong H_{2}(V) \cong \mathbb{Z}$.
Let $N = M \bigcup_{\partial M} D^{6}$, then $N$ is a $6$-dimensional closed manifold and the classification of $M$ is
equivalent to the classification of $N$. Since $N$ is a simply-connected $6$-dimensional closed manifold with $H_{2}(N)$
free, according to the classification results of \cite{Wall} and \cite{Zh}, the Euler characteristic $\chi(N)$, the Pontrjagin class
$p_{1}(N)$, the Stielfel-Whitney class $w_{2}(N)$, and the cubic form
$$H^{2}(N) \times H^{2}(N) \times H^{2}(N) \to \mathbb{Z}$$ are
complete invariants of $N$. From the above lemmae, it is seen that all these invariants are determined by $d$ and $k$.
\end{proof}

\section{Examples}
In this section hypersurfaces fulfilling the conditions are given. Namely, we construct two families of cubic hypersurfaces in $\cp$ with a unique singularity of type $A_{5}$. Then corollary 1 implies they are homeomorphic.

At first we need a recognition principle
to judge if a given singularity is of type $A_{k}$. The principle developed in \cite{BW} can be generalised to our situation without
essential difficulties. We follow the notions from \cite{Ar}. 

A polynomial $f:(\mathbb{C}^{n}, 0) \to (\mathbb{C}, 0)$ is said to be quasihomogeneous of degree $d$ with weights
$(\alpha_{1},\cdots , \alpha_{n})$ if for any $\lambda >0$ we have
$$f(\lambda ^{\alpha_{1}}x_{1}, \cdots , \lambda ^{\alpha_{n}}x_{n})=\lambda ^{d}f(x_{1}, \cdots , x_{n}).$$
A quasihomogeneous polynomial is said to be non-degenerate if $0$ is an isolated singularity of $f$. For example, the normal form of the 
$A_{k}$-singularity
$$z_{1}^{k+1}+z_{2}^{2}+z_{3}^{2}+z_{4}^{2}$$
is a weighted homogeneous polynomial of weights $(\frac{1}{k+1}, \frac{1}{2}, \frac{1}{2}, \frac{1}{2})$ and  degree $1$.

A polynomial $f:(\mathbb{C}^{n}, 0) \to (\mathbb{C}, 0)$ is said to be semiquasihomogenous of degree $d$ with
weights $(\alpha_{1},\cdots , \alpha_{n})$ if it is of the form $f=f_{0}+g$, where $f_{0}$ is a non-degenerate
quasihomogeneous polynomial of degree $d$ with weights $(\alpha_{1},\cdots , \alpha_{n})$ and $g$ is a polynomial
of degree greater than $d$ according to these weights.

\begin{lem}
If $f(z_{1}, z_{2}, z_{3}, z_{4})$ is a semiquasihomogeneous polynomial of degree $1$ with weights
$(\frac{1}{2m}, \frac{1}{2}, \frac{1}{2}, \frac{1}{2})$, then by a change of co-ordinates we can reduce the degree $1$ part to the normal form of 
$A_{2m-1}$-singularity given above and the resulting function will remain semiquasihomogeneous.
\end{lem}

\begin{proof}
There are 10 monomials of degree $1$ with weights $(\frac{1}{2m}, \frac{1}{2}, \frac{1}{2}, \frac{1}{2})$:
$$z_{1}^{2m},\ z_{1}^{m}z_{2},\ z_{1}^{m}z_{3},\ z_{1}^{m}z_{4},\ z_{2}^{2},\ z_{2}z_{3},\ z_{2}z_{4},\ z_{3}^{2},\
z_{3}z_{4},\ z_{4}^{2}.$$
Thus the degree $1$ part of the polynomial is a linear combination of these monomials:
$$a_{1}z_{1}^{2m}+(a_{2}z_{2}+a_{3}z_{3}+a_{4}z_{4})z_{1}^{m}+(a_{5}z_{2}^{2}+a_{6}z_{2}z_{3}+a_{7}z_{2}z_{4}+a_{8}z_{3}^{2}+a_{9}z_{3}z_{4}+a_{10}z_{4}^{2})$$
If $a_2,a_3, a_4\ne0$, then let $w_{2}=a_{2}z_{2}+a_{3}z_{3}+a_{4}z_{4}$. Since $0$ is an isolated singularity, by a change of co-ordinates, we can reduce the quadratic form of $z_{2}$, $z_{3}$, $z_{4}$ to the normal form $z_{2}^{2}+z_{3}^{2}+z_{4}^{2}$. The resulting polynomial then is
$$\begin{array}{cl}
 & a_{1}z_{1}^{2m}+(a_{2}z_{2}+a_{3}z_{3}+a_{4}z_{4})z_{1}^{m}+z_{2}^{2}+z_{3}^{2}+z_{4}^{2}\\
 & \\
= & (a_1-\frac{a_2^2+a_3^2+a_4^2}{4})z_1^{2m}+ (\frac{a_2}{2}z_1^m+z_2)^2+(\frac{a_3}{2}z_1^m+z_3)^2+(\frac{a_4}{2}z_1^m+z_4)^2 \\
\end{array}$$ 

\noindent Then let $w_i=\frac{a_i}{2}z_1^m+z_i$ for $i=2,3,4$. A final re-scaling will reduce the degree $1$ part to the standard form and one easily checks that the other terms will still have weights $>1$.
\end{proof}

\begin{lem}[recognition principle]\label{normal}
If $f(z_1,\dots,z_4)$ is as in lemma 9, then $f$ is equivalent to the normal form $z_{1}^{k+1}+z_{2}^{2}+z_{3}^{2}+z_{4}^{2}$.
\end{lem}

\begin{proof}
By lemma 9 we can reduce $f$ to the normal form plus terms of degree $>1$. According to \cite[page 194, Theorem]{Ar}, a further change of co-ordinates will 
reduce $f$ to $z_{1}^{k+1}+z_{2}^{2}+z_{3}^{2}+z_{4}^{2}+\sum_{1}^{s}c_ie_i$, where $c_i\in \mathbb{C}$ and $e_1, \dots, e_s$ are superdiagonal elements. 
But the Jacobian ideal of the normal form of $A_k$-singularity has generators $z_1^k, z_2, z_3, z_4$, so $1, z_1,\dots, z_1^{k-1}$ is a monomial basis of the local ring $\mathcal{O}/Jf$. The weights are respectively $0, \frac{1}{k+1},\dots, \frac{k-1}{k+1}$. Therefore there is no superdiagonal element.
\end{proof}

\

We construct families of cubic hypersurfaces in $\mathbb{C}\mathrm{P}^{4}$ with a unique $A_{5}$-singularity from cubic surfaces in
$\mathbb{C}\mathrm{P}^{3}$. There are two constructions.

\begin{lem}
Let $V_{0} \subset \mathbb{C}\mathrm{P}^{3}$ be a cubic surface defined by $F_{0}(x_{0},x_{1},x_{2},x_{3})$ such that $P_{0}=[1,0,0,0]$ is a
unique singularity of $V_{0}$. Suppose that the affine curves 
$F_{0}(0,1,x_{2},x_{3})$, $F_{0}(0,x_{1},1,x_{3})$, $F_{0}(0,x_{1},x_{2},1)$ are all smooth. Let $F=F_{0}+x_{0}x_{4}^{2}$, then $F$ defines a
cubic hypersurface $V \subset \mathbb{C}\mathrm{P}^{4}$ with $P=[1,0,0,0,0]$ a unique singularity of $V$ and $P$ is the stabilisation of $P_{0}$.
\end{lem}

\begin{proof}
On the co-ordinate chart $x_{0}=1$, the affine equation is 
$$F_{0}(1,x_{1},x_{2},x_{3})+x_{4}^{2}=0.$$
It is easy to see that $(0,0,0,0)$ is the unique singularity, which is the stabilisation of $P_{0}$. The smoothness of the curves
$F_{0}(0,1,x_{2},x_{3})$, $F_{0}(0,x_{1},1,x_{3})$ and $F_{0}(0,x_{1},x_{2},1)$ ensures that there are no other singularities.
\end{proof}

Now consider the familiy of hypersurfaces defined by the polynomials
$$F(x_{0}, x_{1}, x_{2}, x_{3}, x_{4})=x_{0}(x_{4}^{2}+x_{1}x_{2})+x_{1}x_{3}(x_{1}+ax_{3})+bx_{2}^{3}, (a\ne 0, b\ne 0).$$
According to the classification of cubic surfaces in \cite{BW}, the polynomials
$$F_{0}=x_{0}x_{1}x_{2}+x_{1}x_{3}(x_{1}+ax_{3})+bx_{3}^{2}, (a\ne 0, b\ne 0)$$
define a familiy of cubic surfaces with a unique singularity at $[1,0,0,0]$ of type $A_{5}$. It is easily checked that $F_{0}$ fulfills the conditions
in the above lemma. Therefore $F$ defines a familiy of hypersurfaces of degree $3$ in  $\mathbb{C}\mathrm{P}^{4}$ with $P=[1,0,0,0,0]$ a unique
singularity of type $A_{5}$.

\begin{lem}
Let $V_{0} \subset \mathbb{C}\mathrm{P}^{3}$ be a cubic surface defined by $F_{0}(x_{0},x_{1},x_{2},x_{3})$ such that $P_{0}=[1,0,0,0]$ is a 
unique singularity of $V_{0}$ of type $A_{1}$. Suppose that the affine curve 
$F_{0}(1,0,x_{2},x_{3})$ has a unique singularity at $(0,0)$ and that the affine curves $F_{0}(x_{0},0,1,x_{3})$, $F_{0}(x_{0},0,x_{2},1)$ are
 smooth. Let $F=F_{0}+x_{1}x_{4}^{2}$, then $F$ defines a 
cubic hypersurface $V \subset \mathbb{C}\mathrm{P}^{4}$ with $P=[1,0,0,0,0]$ a unique singularity of type $A_{5}$.
\end{lem}

\begin{proof}
According to \cite{BW}, $F_0$ is of the form $F_0=x_0(x_2^2-x_1x_3)+f_3(x_0,\dots,x_3)$ where $f_3$ is of degree $3$. Let $x_3'=x_3-x_4^2$, then $F(1,x_1,x_2,x_3',x_4)$ is semiquasihomogeneous of $x_1,x_2,x_3',x_4$ of weights $(\frac{1}{2},\frac{1}{2},\frac{1}{2},\frac{1}{6})$. Then according to the recognition principle, $[1,0,0,0,0]$ is an $A_5$-singularity.
\end{proof}

Consider the polynomials
$$F_{0}=x_{0}(x_{2}^{2}-x_{1}x_{3})+ax_{1}^{3}+bx_{3}^{3}, (a\ne 0, b\ne 0).$$
According to \cite{BW}, these polynomials define a familiy of cubic surfaces with a unique singularity at $[1,0,0,0]$ of type $A_{1}$.
It is easily checked that $F_{0}$ fulfills the conditions in the above lemma. Therefore the polynomials
$$F=x_{0}(x_{2}^{2}-x_{1}x_{3})+ax_{1}^{3}+bx_{3}^{3}+x_{1}x_{4}^{2}$$
define a familiy of hypersurfaces of degree $3$ in
$\mathbb{C}\mathrm{P}^{4}$ with $P=[1,0,0,0,0]$ a unique singularity of type $A_{5}$.

It is not clear to the author whether these two families can be deformed to each other.
However, corollary 1 implies that all these hypersurfaces are homeomorphic, and the homeomorphisms are
diffeomorphisms on the nonsingular parts.

Mathematisches Institut, Universtit\"at Heidelberg

Im Neuenheimer Feld 288, D-69120, Heidelberg, Germany

suyang@mathi.uni-heidelberg.de

\end{document}